\documentclass [11pt] {article}

\usepackage{amsmath,amssymb}

\newtheorem{theorem}{Theorem}
  \newtheorem{lemma}{Lemma}

\def\Z{{\mathbb{Z}}}

\def\Q{{\mathbb{Q}}}
\def\K{{\mathbb{K}}}
\def\be{\begin{eqnarray*}}
\def\ee{\end{eqnarray*}} 
\def\bee{\begin{eqnarray}}
\def\eee{\end{eqnarray}}

\begin{document}

\title{{\Large{\bf On the torsion group of elliptic curves induced by Diophantine triples over quadratic fields}}}
{\let\thefootnote\relax\footnotetext{Authors were supported by the Croatian Science Foundation under the project no. 6422.
A.D. acknowledges support from the QuantiXLie Center of Excellence.}}

\author{Andrej Dujella, Mirela Juki\'c Bokun and Ivan Soldo}

\date{}
\maketitle

\begin{abstract} The possible torsion groups of elliptic curves induced by
Diophantine triples over quadratic fields, which do not appear over $\Q$, are
$\Z/2\Z \times \Z/10\Z$, $\Z/2\Z \times \Z/12\Z$ and $\Z/4\Z \times \Z/4\Z$.
In this paper, we show that all these torsion groups indeed appear over some
quadratic field. Moreover, we prove that there are infinitely many Diophantine triples
over quadratic fields which induce elliptic curves with these torsion groups.
\end{abstract}

\medskip
{\small

}

\section{Introduction}
\setcounter{equation}{0}

A set $\{a_1, a_2, \dots, a_m\}$ of $m$ distinct nonzero rationals
is called {\it a rational Diophantine m-tuple} if $a_i a_j+1$ is a perfect square
for all $1\leq i < j\leq m$.
The Greek mathematician Diophantus of Alexandria was the first who studied the existence of
Diophantine quadruples in rationals. He discovered a rational Diophantine quadruple
$\{\frac{1}{16}, \frac{33}{16}, \frac{17}{4}, \frac{105}{16} \}$.
The first example of a Diophantine quadruple in integers, the set $\{1, 3, 8, 120\}$,
was found by Fermat. In 1969, Baker and Davenport \cite{B-D} proved that Fermat's set
cannot be extended to a Diophantine quintuple in integers. It was proved in \cite{D-crelle} that there
does not exist a Diophantine sextuple of integers and there are only finitely many
such quintuples. Euler extended the Fermat set to a rational quintuple $\{1, 3, 8, 120, \frac{777480}{8288641}\}$, and found infinitely
many rational Diophantine quintuples (see \cite{Hea}). The first example of rational Diophantine
sextuple was given by Gibbs in \cite{Gibbs}. It was the set
$\{\frac{11}{192}, \frac{35}{192}, \frac{155}{27}, \frac{512}{27}, \frac{1235}{48}, \frac{180873}{16}\}$.
Recently, in \cite{DujeKaza}, it was proved that there are infinitely many rational Diophantine
sextuples.

 The problem of extendibility and existence of Diophantine $m$-tuples is closely connected with the properties of elliptic curves associated with them.
Let $\{a,b,c\}$ be a rational Diophantine triple. This means that there exist nonnegative rationals
$r,s,t$ such that
\bee\label{eq1}
ab+1=r^2, \quad ac+1=s^2, \quad bc+1=t^2.
\eee

If we want to extend Diophantine triple $\{a,b,c\}$ to a quadruple, we have to solve the system
\bee\label{eq2}
ax+1= \square, \quad bx+1=\square, \quad cx+1=\square.
\eee
The elliptic curve assigned to the system \eqref{eq2} is
\bee\label{e1}
E:\qquad  y^2=(ax+1)(bx+1)(cx+1).
\eee
We say that the elliptic curve $E$ is induced by Diophantine triple $\{a,b,c\}$.

By Mazur's theorem \cite{Mazur}, there are at most four possibilities for the torsion group
over $\Q$ for such curves: $\Z/2\Z \times \Z/2\Z$, $\Z/2\Z \times \Z/4\Z$, $\Z/2\Z \times \Z/6\Z$
and $\Z/2\Z \times \Z/8\Z$. In \cite{D-MW}, it was shown that all these torsion groups actually appears. Moreover, it was shown that every elliptic curve with torsion
group $\Z/2\Z \times \Z/8\Z$ is induced by Diophantine triple (see also \cite{CG}). Questions about the ranks of elliptic curves induced by Diophantine triples was
studied in several articles (\cite{ADP, D2000,D-MW, DP2014}). In particular, such curves were used
for finding elliptic curves with the largest known rank with torsion group
$\Z/2\Z \times\Z/4\Z$ (\cite{DP2014}).

In this paper, we study possible torsion groups of elliptic curves induced by
Diophantine triples over quadratic fields, i.e., we will suppose that $a,b,c$ are elements of
some quadratic field and $ab+1$, $ac+1$ and $bc+1$ are squares in the field.
According to  \cite{Kam,KeMo}, possible torsion groups for such curves which do not appear over $\Q$ are  $\Z/2\Z \times \Z/10\Z$, $\Z/2\Z \times \Z/12\Z$ and $\Z/4\Z \times \Z/4\Z$.
The last one can only appear over quadratic field $\Q(i)$. We show that all these torsion groups indeed appear over some quadratic field and that there are infinitely many
Diophantine triples over quadratic fields which induce elliptic curves with these torsion groups.
In Section \ref{sec:3}, for each of three considered torsion groups,
we present particular examples of elliptic curves over quadratic fields with reasonably large rank.

\section{Torsion groups}

The coordinate transform $x\mapsto\displaystyle \frac{x}{abc}, y\mapsto\displaystyle\frac{y}{abc}$ applied on the curve $E$
leads to the elliptic curve
\bee\label{e2}
E':\qquad  y^2&=&(x + ab)(x + bc)(x + ac)\\
&=& x^3 + (ab+bc+ac)x^2 + (a^2bc + ab^2c + abc^2)x + a^2b^2c^2 \nonumber
\eee
in the Weierstrass form. There are three rational points on $E'$ of order 2:
\[
T_1 = [-ab, 0], \quad T_2 = [-bc, 0],\quad  T_3 = [-ac, 0],
\]
and other two obvious rational points
\[
P = [0, abc],\quad Q = [1, rst].
\]

\subsection{Torsion group $\Z/2\Z \times \Z/10\Z$} \label{sec:2.1}

A method for construction of Diophantine triples was known already to Euler, and in
details is described in \cite{DujeGM}. From that  method and  \eqref{eq1} we have
\bee \label{eq3}
b=\frac{r^2-1}{a}, \quad c=a+b+2r,\qquad r\not\in\{-1,1, 1-a,-1-a\}.
\eee

Our intention is to construct an elliptic curve $E'$ on which the point $P=[0, abc]$ will be of order 5,
i.e., $5P=\mathcal{O}$.

\begin{lemma} \label{Lemma 1}
For the point $P=[0,abc]$ on $E'$ satisfying condition \eqref{eq3}
it holds $5P=\mathcal{O}$ if and only if
\be\label{eq4}
&&(-4 r^2 + 4 r^4)a^4
+ (4 r - 20 r^3 + 16 r^5) a^3
+ ( -1 + 16 r^2 - 40 r^4 + 24 r^6)a^2\\\nonumber
&&\,\,+ (-4 r + 24 r^3 - 36 r^5 + 16 r^7)a
-4 r^2 + 12 r^4 - 12 r^6 + 4 r^8=0.
\ee
\end{lemma}

\noindent
{\bf Proof}: The statement of Lemma~\ref{Lemma 1} follows directly from the condition $\psi_5(P)=0$, where
\[
mP=\left(\frac{\phi_m(P)}{\psi_m(P)^2}, \frac{\omega_m(m)}{\psi_m(P)^3}\right),
\]
and $\phi_m,\psi_m, \omega_m$ are multiplication polynomials \cite[Section 1.3]{SZ}.
\hfill $\Box$

\begin{theorem}\label{Theorem1}
There exist infinitely many quadratic fields $\K$ such that for each of them there exist infinitely many
Diophantine triples which induce elliptic curves with torsion $\Z/2\Z \times \Z/10\Z$ over $\K$.
\end{theorem}

\noindent
{\bf Proof}: If $r^2-1$ is a perfect square, then the polynomial from Lemma~\ref{Lemma 1} factorizes as the product of two quadratic factors in $a$. 
This condition is satisfied for $r=\displaystyle\frac{t^2+1}{2t}$, $t\in \Q \setminus \{ 0 \}$, 
which gives 
\begin{eqnarray*}
& ((4t^6-4t^2)a^2+(4t^7+4t^5+4t^3-4t)a+t^8-2t^6+2t^2-1) \\
& \times ((4t^6-4t^2)a^2+(4t^7-4t^5-4t^3-4t)a+t^8-2t^6+2t^2-1).
\end{eqnarray*}
The roots of the factors are elements of quadratic fields $\Q(\sqrt{t^4+t^2-1})$, respectively $\Q(\sqrt{-t^4+t^2+1})$.
Note that $t^4+t^2-1$ and $-t^4+t^2+1$ are not rational squares, since both conditions $t^4+t^2-1=\tilde{y}^2$ and $-t^4+t^2+1=\tilde{y}^2$ are birationally equivalent to the elliptic curve
\[
E_1:\qquad y^2= x^3 + x^2 + 4 x + 4,
\]
with the torsion group $\Z/6\Z$ and the rank equal to 0 over $\Q$ (and torsion points correspond
to $t=\pm 1$ for the first curve and $t=0,\pm 1$ for the second curve, 
which are excluded values since $t=\pm 1$ gives $r=\pm 1$).

In what follows, let $u \in \Q \setminus \{-1,0,1\}$. The elliptic curve $E_1$ over the quadratic field $\Q(v)$, $v=\sqrt{(2u^2+2u+1)(1-u^2)}$, contains the point
\[
P_1=\left[\frac{-6u-4}{u-1}, \frac{10v}{u-1}\right].
\]
We will show that for all but finitely many rational numbers $u$ the point $P_1$
is a point of infinite order on the curve $E_1$ over $\Q(v)$.
Let us consider under which conditions the point $P_1$ will be a torsion point on
the elliptic curve $E_1$ over $\Q(v)$.
To check that, we try to find quadratic fields
$\mathbb{K}$ for which $E_1(\mathbb{K})_{\rm tors}$ is strictly larger than
$E_1(\mathbb{Q})_{\rm tors}$.

According to \cite{Gonz}, for an elliptic curve with the torsion group $\Z/6\Z$
over $\Q$, its torsion group over a quadratic field can be extended to $\Z/12\Z$,
over at most two quadratic fields,
and to $\Z/3\Z \times \Z/6\Z$, $\Z/2\Z \times \Z/6\Z$ over at most one quadratic field, each.
\begin{itemize}
\item We can use \cite[Proposition 12]{Gonz} (see also \cite{Jeon}) to check if
torsion group $\Z/6\Z$ extends to torsion group $\Z/12\Z$ over some quadratic field.
Transformation $x\mapsto x+4$, $y\mapsto y$, transform  the elliptic curve $E_1$ to
$$E_1':\qquad y^2=x^3+13x^2+60x+100,$$
where the point $[0,10]$ is point of order 6. From the mentioned proposition we conclude that there does not exist a quadratic field over which $E_1$ has torsion $\Z/12\Z$.

\item According to \cite{KeMo}, the torsion $\Z/3\Z \times \Z/6\Z$ can appear only over $\Q(\sqrt{-3})$.
It is easy to show that over $\Q(\sqrt{-3})$ elliptic curve $E_1$ has torsion $\Z/6\Z$ and rank 0.

\item Over $\Q(i)$, the curve $E_1$ has torsion group $\Z/2\Z \times \Z/6\Z$ and rank 0. The point $P_1$
is a torsion point over $\Q(i)$ for
\[
u \in \left\{ -\frac{2}{3}, -1-i, -1+i, -\frac{1}{2}(1-i), -\frac{1}{2}(1+i)\right\},
\]
and in all other cases $P_1$ is a point of infinite order over $\Q(i)$.
\end{itemize}

We conclude that for every rational parameter $u \notin \{-1,-\frac{2}{3},0,1\}$  the point $P_1$ is a
point of infinite order on the curve $E_1$ over $\Q(v)$.
Every multiple $mP_1$ of $P_1$ generates a Diophantine triple such that the
induced elliptic curve $E$ over $\Q(v)$ has torsion $\Z/2\Z \times \Z/10\Z$.
Therefore, there exist
infinitely many Diophantine triples which induce elliptic curves
with torsion $\Z/2\Z \times \Z/10\Z$, over $\Q(v)$.

To show that there are infinitely many different quadratic fields obtained with this construction,
note that, according to Siegel's theorem \cite[Chapter IX.3]{AEC}, there are only finitely many
integer points on the curve
\[
dy^2=(2u^2+2u+1)(1-u^2),
\]
for fixed square-free integer $d$.

\hfill $\Box$


\subsection{Torsion group $\Z/2\Z \times \Z/12\Z$} \label{sec:2.2}

We prove the following theorem:

\begin{theorem}\label{Theorem2}
There exist infinitely many quadratic fields $\K$ such that for each of them there exist infinitely many
Diophantine triples which induce elliptic curves with torsion $\Z/2\Z \times \Z/12\Z$ over $\K$.
\end{theorem}

\noindent
{\bf Proof}:
By the results presented in \cite{DujeKaza}, Diophantine triples $\{a,b,c\}$, where
\bee\label{d3}
a &=& \frac{18t(t^2 - 1)}{(t^2 - 6t + 1)(t^2 + 6t + 1)},\nonumber\\
b &=& \frac{(t - 1)(t^2 + 6t + 1)^2}{6t(t + 1)(t^2 - 6 t + 1)},\\
c &=& \frac{(t + 1) (t^2 - 6 t + 1)^2}{6 t (t - 1) (t^2 + 6 t + 1)},\nonumber
\eee
induce elliptic curves with torsion $\Z/2\Z \times \Z/6\Z$ and positive rank (the point of infinite order is $P=[0,abc]$), for rationals $t\neq -1,0,1$. In order to obtain a point
of the order $12$ on that curve, we consider the point of the order $6$ which has $x$-coordinate
\bee\label{eq5a}
x (P_6) =\frac{2 t^4+ 3 t^3 - 14 t^2+ 3 t+ 2}{3 t (t^2 - 6 t +1)}.
\eee
According to 2-descent proposition \cite[Proposition 4.2.]{Knapp}, we obtain conditions
\bee\label{eq6}
(t^2- 6 t +1)(t^2+ 18 t +1)&=&u^2,\\
6t(t^2+1)&=&v^2.\label{eq7}
\eee
The condition \eqref{eq6} leads to the elliptic curve
\bee\label{e3}
y^2=x^3-x^2-225x-1215,
\eee
which has torsion group $\Z/2\Z$ and rank equal to 1 over $\Q$. The point
of infinite order is $P_{\infty}=[27, -108]$.

For every parameter $t$ which is generated by multiples $mP_{\infty}$ ($m\geq 2$)
of the point $P_{\infty}$,
from \eqref{eq7} we get the elliptic curve over the quadratic field $\Q(\sqrt{6t (1+t^2)})$. Therefore,
over this quadratic field the Diophantine triple \eqref{d3} induces the elliptic curve with torsion $\Z/2\Z \times \Z/12\Z$.

 We conclude that there are infinitely many such quadratic fields,
 since, according to Falting's theorem \cite{Fal}, for fixed square-free integer $d$
 there are only finitely many rational points on the curve
 \[
dy^2=6t(t^2+1)(t^2- 6 t +1)(t^2+ 18 t +1),
\]
of genus 3.

This completes the proof of the Theorem~\ref{Theorem2}.
\hfill $\Box$


\subsection{Torsion group $\Z/4\Z \times \Z/4\Z$} \label{sec:2.3}

We show that there are infinitely many Diophantine triples which induce elliptic curves with torsion
$\Z/4\Z \times \Z/4\Z$ and we give their parameterization.

\begin{theorem}\label{Theorem3}
There exist infinitely many Diophantine triples which induce elliptic curves with torsion
$\Z/4\Z \times \Z/4\Z$ over $\Q(i)$.
Moreover, there exists a Diophantine triple which induces an elliptic curve with torsion
$\Z/4\Z \times \Z/4\Z$ and positive rank over $\Q(i)(t)$.
\end{theorem}

\noindent
{\bf Proof}:
Diophantine triples $\{a,b,c\}$, where
\bee\label{eq11}
a &=& \frac{tu+1}{t-u},\nonumber\\
b &=& -\frac{1}{a},\\
c &=& \frac{4tu}{(tu+1)(t-u)},\nonumber
\eee
induce elliptic curves with torsion $\Z/2\Z \times \Z/4\Z$, over $\Q$, for admissible rational $t,u$ (see \cite{D-MW}).
As in proof of the Theorem~\ref{Theorem2}, we apply 2-descent proposition to get torsion
$\Z/4\Z \times \Z/4\Z$. We obtain conditions
\bee\label{eq8}
-(tu-1)^2&=&m^2,\\
(t^3+t)u^3+(t^3+t)u&=&n^2.\label{eq9}
\eee
The condition \eqref{eq8} is always satisfied over $\Q(i)$. From \eqref{eq9}, by using
substitutions $U=(t^3+t)u, N=(t^3+t)n$, we get the elliptic curve
\[
U^3+(t^3+t)^2 U=N^2,
\]
with the point $P_t=[t^2+1, (t^2+1)^2]$. It is easy to check that the point $P_t$
does not generate an appropriate Diophantine triple, but its multiples $mP_t, m>1$, do.
For $m=2$, we obtain
\bee\label{eq10}
u=\frac{(t^2-1)^2}{4 t (t^2+1)}.
\eee
Inserting \eqref{eq10} into \eqref{eq11},
we obtain a parametric family of Diophantine triples $\{a,b,c\}$ where
\be
a&=&\frac{t(t^4+ 2 t^2+5)}{3 t^4 + 6 t^2-1},\\
b&=&\frac{- 3 t^4- 6 t^2+1}{t^5+ 2 t^3+5t},\\
c&=&\frac{16 t (t^4-1)(t^2-1)}{(t^4+ 2 t^2+5) (3 t^4+ 6 t^2-1)},
\ee
which induces an elliptic curve with torsion $\Z/4\Z \times \Z/4\Z$,
over $\Q(i)(t)$. By applying Silverman's specialization theorem \cite[Theorem 11.4]{ATAEC} (using, for example, the specialization $t=2$) it can be
shown that the point  $P=[0,abc]$ is point of infinite order.
\hfill $\Box$


\section{Ranks} \label{sec:3}

By the results of Kenku and Momose \cite{KeMo} and Kamienny \cite{Kam},
there are exactly 26 possibilities for the torsion group of elliptic curves
defined over quadratic fields. Several authors constructed elliptic curves
with reasonably large rank and given torsion group over quadratic fields
(see \cite{ADJBP,BBDN,DJB,JB,Naj-Ov,Rab}). The current records can be found on the
web page \cite{D-torsquad}. Here we show that in the cases of three torsion groups
considered in the previous section, curves with positive rank can be induced
by Diophantine triples over certain quadratic fields. In computations of the ranks,
we will use {\tt mwrank} \cite{mwrank} and {\tt Magma} \cite{magma}.
Since our curve will have rational coefficients, and in order to determine
their rank over a quadratic field $\Q(\sqrt{d})$ we will use the formula
(see \cite{Birch})
$$ {\rm rank}(E(\Q(\sqrt{d})) = {\rm rank}(E(\Q)) + {\rm rank}(E^{(d)}(\Q)), $$
where $E^{(d)}$ denotes the $d$-quadratic twist of $E$.

\subsection{Torsion group $\Z/2\Z \times \Z/10\Z$}

By taking $u=3$, i.e. $v=-\frac{25}{8}$, in construction from Section \ref{sec:2.1}, we
get the curve with rank $3$ and torsion group $\Z/2\Z \times \Z/10\Z$ over
$\Q(\sqrt{-2})$. In this case, $P_1=[-11,25\sqrt{-2}]$, $t=\frac{2}{5}\sqrt{-2}$
and $a=\frac{475}{561}+\frac{12737}{22440} \, \sqrt{-2}$.
Hence, the curve is induced by the Diophantine triple
$$ \left\{ \frac{475}{561}+\frac{12737}{22440} \, \sqrt{-2},
-\frac{475}{561}+\frac{12737}{22440} \,\sqrt{-2}, \frac{160}{561} \sqrt{-2} \right\}. $$
The Weierstrass equation of the curve is
{\small
$$ y^2 = x^3+x^2-61404142096090881x-20861928799251086002759425. $$ }%
Three independent points of infinite order are:
\begin{eqnarray*}
& [865303425, 23956226997120], \\
& \displaystyle{\left[\frac{48954515537984337}{16008001}, \frac{10791931818384647817975000}{64048012001} \right]}, \\
& \displaystyle{\left[\frac{86963667871383}{299209}, \frac{435438077091034960800}{163667323} \, \sqrt{-2}\right]}.
\end{eqnarray*}
Let us mention that the current record for the
rank for arbitrary elliptic curves over quadratic fields with torsion group
$\Z/2\Z \times \Z/10\Z$ is $4$ (see \cite{BBDN,ADJBP}).

\subsection{Torsion group $\Z/2\Z \times \Z/12\Z$}

By taking multiples $2P_{\infty}$ and $3P_{\infty}$ from Section \ref{sec:2.2},
we get $t=6/35$ and $t=41615/426$, and in both cases we obtain curves with rank
between $1$ and $3$. For larger multiples, the coefficients of the curves are too large
and we are not able to compute reasonable bounds for the rank.
For $t=41615/426$, we can conclude that the rank of the corresponding curve
over $\Q(\sqrt{5117449349905165})$ is equal to $3$ assuming the Parity conjecture.

We can use other parametric families of rational Diophantine triples
which satisfy the condition of \cite[Lemma 1]{DujeKaza}, and hence
induce elliptic curves over $\Q$ with torsion $\Z/2\Z \times \Z/6\Z$,
as the starting point for our construction, e.g.,
\[
a=\frac{u^3-9u}{6(u^2-1)}, b=-\frac{9(u^2-1)}{2(u^3-9u)}, c=-\frac{16u(u^2-3)}{3(u^4-10u^2+9)}.
\]
 Now the condition for the torsion $\Z/2\Z \times \Z/12\Z$ becomes
 $3(u-1)(u+1)(u^2+15) = v^2$. Thus, we may take here $u=-7$, and we get the curve
 induced by the Diophantine triple
 $$  \left\{ -\frac{35}{36}, \frac{27}{35}, \frac{161}{180} \right\},$$
which has torsion group $\Z/2\Z \times \Z/12\Z$ and rank (unconditionally) equal to $2$
over $\Q(\sqrt{-155})$.  The Weierstrass equation of the curve is
$$ y^2+xy+y = x^3-49428958x+130902669056, $$
and two independent points of infinite order are:
$$ [-2510, -487783], \displaystyle{\left[-\frac{95078581}{245},
\frac{166483532709}{8575} \, \sqrt{-155}\right]}. $$
The current record for the
rank for arbitrary elliptic curves over quadratic fields with torsion group
$\Z/2\Z \times \Z/12\Z$ is $4$ (see \cite{BBDN}).

\subsection{Torsion group $\Z/4\Z \times \Z/4\Z$}

Inserting $t=4/3$ in the family of Diophantine triples from Section \ref{sec:2.3},
we obtain the curve with rank $6$ and torsion group $\Z/4\Z \times \Z/4\Z$ over $\Q(i)$.
The curve is induced by the Diophantine triple
$$ \left\{ \frac{3796}{4653}, -\frac{4653}{3796}, \frac{78400}{490633} \right\}. $$
The Weierstrass equation of the curve is
{\small
$$ y^2 = x^3+x^2-1588627573982287131943200x-507161545884329501301628000492040652. $$ }%
Six independent points of infinite order are:
\begin{eqnarray*}
& [-890497354044, 448726623142928130], \\
& [-899563900533, 440419889828558640], \\
& \displaystyle{\left[\frac{2502824381840097811}{632025}, \frac{3736538268665587610111875016}{502459875}\right]}, \\
& [-1089076885194, 262231774368503940 \, i], \\
& [-1926913622169, 2144909371334503410 \, i], \\
& \displaystyle{\left[-\frac{10573435624608518034}{6175225}, \frac{25709440364558354804130497052}{15345434125} \, i\right]}.
\end{eqnarray*}
The current record for the
rank over $\Q(i)$ for arbitrary elliptic curves with torsion group
$\Z/4\Z \times \Z/4\Z$ is $7$ (see \cite{DJB}).

\bigskip
\noindent
{\sc
Andrej Dujella\\
Department of Mathematics\\
University of Zagreb\\
Bijeni\v{c}ka cesta 30\\
10000 Zagreb, Croatia\\
e-mail: {\rm duje@math.hr}

\vspace*{1cm}
\noindent
Mirela Juki\'c Bokun\\
Department of Mathematics\\
University of Osijek\\
Trg Ljudevita Gaja 6\\
31000 Osijek, Croatia\\
e-mail: {\rm mirela@mathos.hr}

\vspace*{1cm}
\noindent
Ivan Soldo\\
Department of Mathematics\\
University of Osijek\\
Trg Ljudevita Gaja 6\\
31000 Osijek, Croatia\\
e-mail: {\rm isoldo@mathos.hr}
}
\end{document}